\documentclass[journal]{IEEEtran}
\usepackage{graphicx} % Required for inserting images
\usepackage{amsmath}
\usepackage{amsfonts}
\usepackage[a4paper, left=2cm, right=2cm, top=2cm, bottom=2cm]{geometry}
\usepackage{algorithm}
\usepackage{algpseudocode}
\usepackage{url}

\usepackage{comment}

\title{Community Detection using the Edge Laplacian Dynamics vis-a-vis Ricci Flow}

\author{
\IEEEauthorblockN{Pratibha Bhandari, Soumyendu Raha}

\IEEEauthorblockA{
Department of Computational and Data Sciences\\
Indian Institute of Science (IISc) Bengaluru, India\\
Email: pratibhab@iisc.ac.in, raha@iisc.ac.in
}
}

\begin{document}
\maketitle
\begin{center}
\textbf{Abstract}
\end{center}

It has been found that utilizing the geometric properties of the graph dynamics can bring out crucial information of the data that a statistical analysis will not. The properties of Ricci flow bring out hidden dynamics of the data in the same way as decomposition of smooth manifolds. In this paper we explore an alternative to Ricci flow, where we find that we can bring out similar properties like it using a simpler and well defined Edge Laplacian to practice community detection in graph. Moreover, the results obtained by the present approach show that a flow driven by the Edge Laplacian yields results similar to those of the Ricci flow. In contrast, computing the Forman-Ricci curvature requires looping over the adjacent edges, while computing Olivier-Ricci curvature involves the computation of Wasserstein distance and probability distribution associated with the graph nodes, making the Edge Laplacian approach computationally more efficient.

\section{Introduction}
Graph representation learning is one of the dominant routes that is kept in use to work on graphs. It utilizes the properties of nodes and edges to perform link prediction, node classification, and Graph classification using different variants of Graph Neural Networks \cite{scarselli2009gnn, velickovic2018gat, kipf2017gcn}. These theories work best on data exhibiting Euclidean or grid-like setup, particularly when the network architecture incorporates the invariance of those setup.  But knowing the geometry behind the characterization of graph can be more beneficial as handled using non-Euclidean domains and manifolds. Using such atrocities can help us understand more about graphs. It is one of the emerging fields in sensor networks, networks in brain imaging, classification of social networks, disease spread networks, regulatory networks in genetics, etc.

Geometric deep learning structures are of major interest for the networks which do not have a grid like structure or abide by the Euclidean dynamics. For data that do not abide the Euclidean dynamics or have grid like structure, curvature based learning facilitates the objective. Most of these concepts are taken from the manifold theory. Hence the whole idea is to apply curvature dynamics to the learning on non-Euclidean data, which takes Riemannian manifolds into picture . 
As we broadly know that a very known concept, CNN(Convolutional Neural Network) originated from the paper \cite{lecun1989backprop}, is applied to grid like structure data in order to use the convolutional layers and downsampling layers. We want to translate this CNN framework to the non-Euclidean or geometric setup for the clustering of graphs. For this, Ricci flow comes into existence.

Hamilton \cite{Hamilton1982} first introduced  Ricci flow in the smooth Riemannian setting; Chow and Luo \cite{ChowLuo2003} later developed a combinatorial/discrete version for triangulated surfaces. Furthermore, Perelman’s contribution ~\cite{Perelman2002,Perelman2003a,Perelman2003b} to Ricci flow theory was transformational, he made it powerful enough to classify 3-manifolds. The new notion of Ricci curvature for metric spaces was introduced by Ollivier in \cite{Ollivier2009}, where he introduced the Ollivier-Ricci curvature based on optimal transportation theory. In this, they considered the network as a counterpart of discrete manifold. This concept has then been used in many fields; Financial market fragility \cite{Sandhu2016}, Gene networks for detection of cancer \cite{Sandhu2015CancerNetworks} (they found that the curvature of the gene network reliably distinguishes cancer from normal samples, with cancer networks showing higher curvature.), Wireless network capacity and network congestion \cite{Wang2014WirelessRicci} \cite{Ni2015RicciInternet} (It explores the fact that negative curvature is considered to be the cause of large queue occupancy, region of limited capacity in a wireless network, hence determining the network backbone and congestion), and many more. Further \cite{Weber2016FormanRicci} came up with Forman Ricci curvature, which is based on the Laplacian, but it is less geometric. Forman Ricci curvature is also known to be the discrete version of Ricci curvature; it solely depends on combinatorial data (adjacency matrix) and prescribed weights of the network. It does not depend on the way the edges of a graph are drawn or the distribution of edge weights. It was also used for noisy graph matching \cite{Ni2018NetworkAlignmentRicci}. Further, \cite{Weber2016FormanRicci} uses the Forman Ricci flow for anomaly detection.

Discrete flow in networks\cite{ni2019community} was motivated by the Hamilton's Ricci flow, used for community detection defined on weighted graphs. It introduced the concept that heavily traveled edges will be stretched and sparsely traveled edges will shrink. This idea was inspired by \cite{JostLiu2014Ollivier}, \cite{Ni2015RicciInternet}, \cite{Samal2018ComparativeRicci} and \cite{Sreejith2017Systematic}.

One of the most famous concepts that takes smooth manifold learning into picture is Ricci Flow for community detection\cite{ni2019community}. Some methods that were widely used for clustering (community detection) traditionally were label propagation \cite{Raghavan2007} (a fast, randomized algorithm in which each node adopts the most frequent label among its neighbors and communities emerge naturally through repeated updates), Random walks \cite{Rosvall2008Maps} (Uses random walks to model information flow in a network, communities are identified as regions where a random walker is likely to stay for longer periods), Edge betweenness Optimization \cite{GirvanNewman2002} (Identifies communities by successively eliminating high-betweenness edges, which exposes the network’s hierarchical structure), Modularity \cite{NewmanGirvan2004} (Identifies communities by maximizing modularity, a measure that compares observed intra-community connections with those expected in a random graph), etc.

 In this work we explore obtaining the same results that a Ricci flow approach produces  but using simpler and natural approach to diffusion and curvature for community detection in graphs. We want to determine the curvature of a graph but not using the complicated Ricci curvature with Optimal Transportation theory. Instead, we utilize a well defiened Edge Laplacian considering heat equation on networks, and see how the edge weight behaves under the perturbation dynamics. Where the dependent variable of heat equation represents the edge perturbation, and hence called perturbation dynamics. And finally removing the perturbation from the weights and analyzing how the curvature (edge weights here) dynamics behave. Also, plotting the major inferences about the graph and the changes in weight dynamics due to the perturbation.

 This paper is structured as follows: section II provides preliminaries providing the whole understanding of manifold theory and heat diffusion implementations, introducing the Ricci flow and its theoretical implementation and finally describing Ricci Forman curvature; section III provides the motivation and contributions in the paper; section IV  introduces the proposed algorithm and its stability; section V explains about the results procured, dataset used and computational efficiency; And finally section VI gives the conclusion of the paper with the involved future work.

\section{Preliminaries}
\subsection{Flow on Manifolds and Graphs}

Many physical, biological, and information-processing systems evolve through \emph{flows}: heat spreads on surfaces, probability mass diffuses over state spaces, and information propagates across networks. When the underlying domain is non-Euclidean such as a curved surface or a graph, classical tools from Euclidean calculus no longer apply directly. Manifolds and graphs provide principled geometric frameworks for modeling such domains, while differential operators define how signals evolve on them.

\subsubsection{Diffusion as the Fundamental Dynamical Process}

We begin with diffusion, as it naturally motivates the operators introduced later. On a smooth domain, heat propagation is governed by the diffusion equation
\begin{equation}
\frac{\partial g(x,t)}{\partial t} = -c\,\Delta g(x,t), 
\qquad g(x,0) = g_0(x),
\label{eq:heat}
\end{equation}
where \( g(x,t) \) denotes the temperature at point \( x \) and time \( t \), the diffusivity constant \( c \), and the Laplace operator \( \Delta \).

The solution admits a spectral representation
\begin{equation}
g(x,t) = e^{-t\Delta} g_0(x)
= \sum_{i \ge 0} \langle g_0, \psi_i \rangle \, e^{-t\lambda_i}\, \psi_i(x),
\end{equation}
where \( \{(\lambda_i, \psi_i)\} \) are the eigen pairs of the Laplacian. This expression shows that diffusion acts as a low-pass filter in the spectral domain, progressively attenuating high-frequency components.

The associated heat kernel is defined as
\begin{equation}
h_t(x,x') = \sum_{l \ge 0} e^{-t\lambda_l}\psi_l(x)\psi_l(x'),
\end{equation}
which measures the amount of heat transferred from \( x \) to \( x' \) over time \( t \). Unlike the Euclidean case, the kernel is not shift-invariant, reflecting the intrinsic geometry of the domain.

From the heat kernel arises the diffusion distance
\begin{equation}
d_t^2(x,x') =
\sum_{l \ge 0} e^{-2t\lambda_i}
\big(\psi_l(x) - \psi_l(x')\big)^2,
\end{equation}
which averages over all possible paths between points and is therefore more robust to structural perturbations than geodesic distance.

\subsubsection{Manifolds as Continuous Geometric Domains}

A manifold is a space that may be globally curved or topologically complex, yet locally resembles Euclidean space. One of a canonical example is Earth; while spherical at large scales, it appears flat in any sufficiently small neighborhood. This local Euclidean structure allows classical calculus to be extended to non-Euclidean domains.

Formally, a smooth  manifold \( M \) (\( d \)-dimensional) is a topological space such that each point \( y \in M \) has a homeomorphic neighbourhood to \( \mathbb{R}^d \). Therefore, associated with each point is a tangent space \( T_yM \), which gives a local linear approximation to the manifold. The set of all such tangent spaces together constitutes the tangent bundle \( T_M \).

Equipping each tangent space with an inner product yields a Riemannian metric, allowing distances, angles, and volumes to be measured intrinsically. Although a manifold may admit embeddings into higher-dimensional Euclidean spaces, such embeddings are generally not unique. Properties determined solely by the Riemannian metric are intrinsic, whereas those depending on a specific
embedding are extrinsic.

\subsubsection{Fields and Differential Operators on Manifolds}

A scalar field on a manifold is a smooth function
\begin{equation}
g : M \rightarrow \mathbb{R},
\end{equation}
while a tangent vector field assigns to each point \( y \in M \) a vector \( G(y) \in T_yM \),
\begin{equation}
G : M \rightarrow T_M.
\end{equation}
Tangent vector fields naturally describe flows of material or information on the manifold.

We define the Hilbert spaces \( \ell_h^2(M) \) and \( \ell_h^2(T_M) \) with inner products
\begin{equation}
\langle g, h \rangle_{\ell_h^2(M)} = \int_M g(y)h(y)\,dy,
\end{equation}
\begin{equation}
\langle G, H \rangle_{\ell_h^2(T_M)} =
\int_M \langle G(y), H(y) \rangle_{T_yM}\,dy.
\end{equation}

Since manifolds do not possess a global vector space structure, differentiation must be defined locally. The differential of a function is a linear operator
\begin{equation}
d_g : T_M \rightarrow \mathbb{R},
\end{equation}
which can be identified with the Riemannian gradient via
\begin{equation}
d_g(y)(G(y)) =
\langle \nabla g(y), G(y) \rangle_{T_yM}.
\end{equation}

The divergence operator
\begin{equation}
\mathrm{div} : \ell_h^2(T_M) \rightarrow \ell_h^2(M)
\end{equation}
is defined as the formal adjoint of the gradient,
\begin{equation}
\langle G, \nabla g \rangle_{\ell_h^2(T_M)}
= \langle -\mathrm{div}\,G, g \rangle_{\ell_h^2(M)}.
\end{equation}

Combining these operators yields the Laplace--Beltrami operator
\begin{equation}
\Delta g = -\mathrm{div}(\nabla g),
\label{eq: div_grad}
\end{equation}
whose associated quadratic form
\begin{equation}
\langle \nabla g, \nabla g \rangle_{\ell_h^2(T_M)}
= \langle g, \Delta g \rangle_{\ell_h^2(M)}
\end{equation}
is known as the Dirichlet energy that measures the smoothness of \( g \) on themanifold.

\subsubsection{Graphs as Discrete Analogues of Manifolds}

Graphs can be viewed as discrete geometric spaces that mirror many properties of manifolds. Consider a weighted undirected graph \( G=(V,E) \), where each vertex \( m,n \in V \) has weight \( a_m>0 \) and each edge \( (m,n)\in E \) has weight \( w_{mn}\ge 0 \). Vertex functions \( f,g:V\rightarrow\mathbb{R} \) and edge functions \( F, G:E\rightarrow\mathbb{R} \) serve as discrete counterparts of scalar and vector fields. The associated inner products are
\begin{equation}
\langle f, g \rangle_{\ell_h^2(V)} = \sum_{m\in V} a_m f_m g_m,
\end{equation}
\begin{equation}
\langle F, G \rangle_{\ell_h^2(E)} =
\sum_{(m,n)\in E} w_{mn} F_{mn} G_{mn}.
\end{equation}

The graph gradient is defined as
\begin{equation}
(\nabla f)_{mn} = f_m - f_n,
\end{equation}
and the graph divergence as
\begin{equation}
(\mathrm{div}\,F)_m =
\frac{1}{a_m} \sum_{n:(m,n)\in E} w_{mn} F_{mn}.
\end{equation}

These operators satisfy an adjointness relation analogous to the continuous case, leading to the graph Laplacian
\begin{equation}
(\Delta f)_m =
\frac{1}{a_m} \sum_{(m,n)\in E} w_{mn}(f_m - f_n).
\end{equation}
In matrix form,
\begin{equation}
\Delta = A^{-1}(D - W),
\end{equation}
where \( A \), \(W\) and \( D \) represent the diagonal matrix of vertex weights, adjacency, and degree matrices, respectively.

\subsection{Discrete Ricci Flow:}
This algorithm implements Ricci curvature utilizing optimal transportation theory to define the probability measure of the vertices connected by an edge.  It is defined on weighted graphs, which keeps on changing for different iterations.  Hence the edges with large curvature in a graph shrinks (representing same community) and edges with small curvature expands (representing different community). Further, removing the stretched edges by considering some threshold for the cut, one can interpret the communities in the evolved structure of the graph.

\subsubsection{Ricci curvature and Optimal transportation} Curvature basically gives us the idea about how space is curved. As stated by Gauss, the curvature depends on how the Riemannian metric is defined on the surface. G. Monge first formulated the optimal transport problem in 1781 \cite{monge1781}. His goal was to minimize the cost of transporting iron ore from a set of mines to a collection of factories. Formally, let $X$ and $Y$ be probability spaces representing the locations of mines and factories, equipped with probability measures $\eta$ and $\zeta$. The cost of moving  a unit of mass from $x \in X$ to $y \in Y$ is given by the function $c(x,y)$. Such that the cost is taken to be the distance $d(x,y)$ when $X = Y$. A transport map $T_r : X \to Y$ pushes forward $\eta$ to $\zeta$, and Monge's formulation seeks a map minimizing the total cost:
\[
\inf_{T_r} \int_X c(x, T_r(x)) \, d\eta(x).
\]

A major breakthrough came in 1930 when Kantorovich reformulated the problem as a linear optimization problem. Instead of transport maps, he considered \emph{transport plans}, namely probability measures $\nu$ on $X \times Y$ with marginals $\eta$ and $\zeta$. The optimal cost becomes
\[
\inf_{\nu \in \Gamma(\eta,\zeta)} 
\int_{X \times Y} c(x,y)\, d\nu(x,y),
\]
where $\Gamma(\eta,\zeta)$ is the set of all admissible transport plans.  Hence, the Wasserstein distance (or Earth Mover’s Distance) 
\(
W(\eta,\zeta)
\) is defined by the optimal transport cost, when $X$ is a metric space and $c(x,y) = d(x,y)$

\bigskip

\subsubsection{Wasserstein distance and Ollivier Ricci curvature}
Wasserstein distance plays a central role in Ollivier’s coarse Ricci curvature. Consider an $n$-dimensional Riemannian manifold $(M^n, d)$ with Riemannian volume measure $\eta$. For $\varepsilon > 0$, at each point $x \in M$ the probability measure $m_x$ is defined as
\[
m_x = \frac{1}{\eta(O(x,\varepsilon))}\,\eta|_{O(x,\varepsilon)},
\]
where $O(x,\varepsilon)$ is the geodesic ball of radius $\varepsilon$ around $x$.

Ollivier showed that for small $\varepsilon$,
\[
W(m_x, m_y) = (1 - \kappa(x,y)) \, d(x,y),
\]
where
\[
\kappa(x,y) = \frac{\mathrm{Ricci(\tau,\tau})}{n-1} \, \varepsilon^2 + o(\varepsilon^2),
\]
and $\tau$ is the tangent vector at $x$ to the geodesic connecting $y$ and $x$. This observation gives rise to a notion of Ricci curvature for general metric measure spaces.

\bigskip

\subsubsection{Ollivier Ricci curvature on general metric spaces.}
Given a metric space $(X,d)$ equipped with a family of probability measures $\{m_x\}_{x\in X}$, the Ollivier Ricci curvature between two points $x,y \in X$ is defined as
\[
\kappa_{xy}
= 1 - \frac{W(m_x, m_y)}{d(x,y)}.
\]

A commonly used family of measures on graphs is the  $(\beta,q)$-lazy random walk distribution with $\beta \in [0,1] $ and $\ q>0$,

\[
m_x^{\beta,q}(x_i) =
\begin{cases}
\beta, 
& \text{if } x_i = x, \\[0.1em]
\dfrac{1-\beta}{b}\, \exp\!\left[- d(x,x_i)^q \right],
& \text{if } x_i \in \pi(x), \\[1em]
0, 
& \text{otherwise},
\end{cases}
\]
where $\pi(x)$ denotes the neighbors of $x$; and $b$ is the normalization constant
\[
b = \sum_{x_i \in \pi(x)} \exp\!\left[- d(x,x_i)^q \right].
\]

\bigskip

\subsubsection{Ricci flow on graphs.}
Given an edge $(x,y)$ with weight $d^{(j)}(x,y)$ at iteration $j$, the discrete Ricci flow updates the weight as
\[
w_{xy}^{(j+1)}
= d^{(j)}(x,y) - \kappa_{xy}^{(j)} \, d^{(j)}(x,y).
\]

\subsection{Forman--Ricci curvature.}
Forman-Ricci curvature, a discretized form of Ricci curvature defined on edges of a weighted graph, introduced by Forman as a combinatorial analogue of the Bochner-Weitzenb\"ock formula.
For an undirected weighted graph $G=(V,E)$ with node weights $w_v$ and edge weights $w_e$, the Forman--Ricci curvature of an edge $e=(u,v)$ is given by
\begin{equation}
\label{eq:forman}
\begin{aligned}
\mathcal{F}(e)
= {} & w_e\Bigg(
\frac{w_u}{w_e} + \frac{w_v}{w_e} - \sum_{e_u \sim e} \frac{w_u}{\sqrt{w_e\,w_{e_u}}}\\
&\quad
- \sum_{e_v \sim e} \frac{w_v}{\sqrt{w_e\,w_{e_v}}}
\Bigg).
\end{aligned}
\end{equation}

where $e_u \sim e$ (resp.\ $e_v \sim e$) denotes edges incident to $u$ (resp.\ $v$), excluding $e$ itself. Positive curvature indicates locally dense connectivity, while negative curvature reflects bottlenecks or bridge-like structures in the graph.

\section{Motivation and Contributions}
\begin{itemize}
    
\item The proposed framework is inspired by the Ricci flow, which exhibits several similarities to the heat equation. Let $g_{ij}$ denote the Riemannian metric on a manifold $M$, and let $R_{ij}$ denote its Ricci curvature tensor. The Ricci flow is governed by the nonlinear second-order partial differential equation

$$\frac{\partial g_{ij}}{\partial t} = -2R_{ij}.$$

Our objective is to derive an analogous perturbation dynamics on graphs of the form

$$\frac{\partial u_e}{\partial t} = -Lu_e$$

where $u_e$ denotes the perturbation associated with the edge $e$ of the graph (or the edge weight-perturbation), and $L$ is the corresponding Edge Laplacian. The proposed formulation is motivated by one of the fundamental properties of Ricci flow, its interpretation as a nonlinear heat equation that progressively smooths the underlying curvature.

\begin{comment}
 The whole idea is inspired from the Ricci flow, the flow having many similarities with the heat equation. let $g_{ij}$ be the Riemannian metric on a manifold $M$, with Ricci curvature as $R_{ij}$, we have a non-linear second order partial differential equation as $$\frac{\partial} { \partial t} g_{ij}=-2R_{ij}$$ which we are trying to convert to a perturbation dynamics as 
 $$\frac{\partial}{\partial t}u_e=-Lu_e$$
 where $u_e$ is the perturbation of $e$ edge present in the graph and $L$ is the well defined edge Laplacian of the graph.  
 Utilising one of the key properties of it that its curvature evaluation according to the non-linear heat equation. 
\end{comment}

\item In Ricci flow, the edge weights are iteratively updated based on the Ricci curvature of the graph. The underlying principle is that edges connecting nodes from different communities tend to increase in weight, whereas edges within the same community tend to decrease in weight during the flow. This behavior naturally strengthens the separation between communities.

Motivated by this curvature-driven separation mechanism, we introduce an edge-perturbation dynamics. Our hypothesis is that the perturbation carries information about the structural role of an edge: positive perturbations correspond to edge contraction and are expected to be more pronounced for intra-community edges, whereas negative or constrained perturbations are expected for inter-community edges. We investigate whether the resulting perturbation dynamics reproduce the community-separating behavior of Ricci flow without explicitly computing Ricci curvature.
\begin{comment}
    Inspired by this idea, we propose a perturbation-based dynamics on graph edges. Our hypothesis is that a larger edge perturbation indicates that the corresponding edge is more likely to connect nodes within the same community, while smaller or constrained perturbations are more likely to occur on edges connecting nodes from different communities.
\end{comment}

The corresponding update rules are

Ricci Flow : $\qquad W_e^{(i+1)} = W_e^{(i)} - \kappa_e W_e^{(i)}$

Proposed Method: \qquad $W_e^{(i+1)} = W_e^{(i)} - u_e^{(i)}$.

where $W_e^{i}$ is the weight of the edge $e$ at $i^{th}$ iteration; $\kappa_e$ is the ricci curvature for the edge $e$.
\begin{comment}
    
\item In Ricci flow, only edge weights are taken into consideration. The theory that it works upon is the edge weight increases for nodes that are not in the same community, and decreases for the intra community settings using the dynamics. Similar to this we utilise the theory that growing edge perturbation implies the intra-community setting and any restriction in the perturbation or small perturbations are only the result of vertices belonging to different communities.

Ricci Flow: $W_e^{i+1}=W_e^i- \kappa_e W_e^i$

Ours: $W_e^{i+1}=W_e^i-u_e$

\end{comment}

\end{itemize}
\section{Proposed Algorithm}

\subsection{Hodge 1-Laplacian}

The Hodge Laplacian provides a natural extension of the graph Laplacian from vertex functions to higher-order simplicial objects. In particular, the Hodge 1-Laplacian acts on edge flows (or 1-cochains) and captures the interaction between edges through their incidence with vertices and triangles. Let $B_1 \in \mathbb{R}^{|V|\times |E|}$ denote the vertex-edge incidence matrix and $B_2 \in \mathbb{R}^{|E|\times |T|}$ denote the edge-triangle incidence matrix, where $V$, $E$, and $T$ represent the sets of vertices, edges, and triangles, respectively. The Hodge 1-Laplacian is defined as
\begin{equation}
    L_1 = L_1^{\mathrm{down}} + L_1^{\mathrm{up}},
\end{equation}
where
\begin{equation}
    L_1^{\mathrm{down}} = B_1^{\top}B_1,
    \qquad
    L_1^{\mathrm{up}} = B_2B_2^{\top}.
\end{equation}
Thus,
\begin{equation}
    L_1 = B_1^{\top}B_1 + B_2B_2^{\top}.
\end{equation}
The down component describes interactions between edges that share a common vertex, whereas the up component describes interactions between edges that belong to a common triangle. Consequently, the Hodge 1-Laplacian provides a topology-aware operator for modeling the dynamics of edge-based quantities and distinguishes gradient, curl, and harmonic components of edge flows.

\subsection{Edge-Based Flow Dynamics}
Beyond node-based diffusion, flows can also be defined directly on graph edges. Motivated by heat diffusion and Ricci flow, we introduce an edge perturbation \( u_e \) evolving according to
\begin{equation}
\frac{du_e}{dt} = -Lu_e,
\label{eq: heat lap}
\end{equation}
where \( L \) is an irrotational Edge Laplacian with zero curl,
\begin{equation}
L = \mathrm{diag}(W_e)^{\frac{-1}{2}} M\,\mathrm{diag}(W_v)\,M^{\top} \mathrm{diag}(W_e)^{\frac{-1}{2}}.
\label{eq:edge_lap}
\end{equation}
Here \( M \) denotes the incidence matrix, while \( W_e \) and \( W_v \) are edge and vertex weight vectors. The above equation (\ref{eq:edge_lap}) is constituted from the two equations introduced in the preliminaries equation (\ref{eq:heat}) and (\ref{eq: div_grad}). Equation (\ref{eq:heat}) comes from the basic definition of diffusion equation, where Laplacian operator (\ref{eq: div_grad}) can be considered as the divergence of the gradient.
\begin{equation}
\frac{\partial u(x,t)}{\partial t} = -c\,\Delta u(x,t),\ \ \  u(x,0) = u_0(x),
\label{eq: perturb}
\end{equation}

$$\Delta u = \mathrm{div}(\nabla u),\ \ \ \ where\ \  \Delta=L$$

$\Delta$ is also called Laplace Beltrami operator for some metric defined on manifold $M$.  The Beltrami Laplacian governs diffusion on vertex-based functions, whereas the Hodge $1$-Laplacian extends this notion to edge-based flows, capturing their interactions through shared vertices and higher-order simplices. Hence, for our algorithm, we use the normalized lower Hodge 1-laplacian considering there are no higher dimensional simplices. After discretizing equation (\ref{eq: perturb}), the edge weights are updated iteratively as
\begin{equation}
W_e^{(i+1)} = W_e^{(i)} - u_e^{(i)},
\end{equation}
allowing the analysis of how local perturbations propagate and reshape the global graph structure over time. Algo[\ref{propAlg}] gives the proposed algorithm replacing the Ricci Curvature with the Edge Laplacian as the diffusion operator.
\begin{figure}[h!]
\hspace{-1.5cm}
   \centering
    \includegraphics[width=9cm , height=3.0cm]{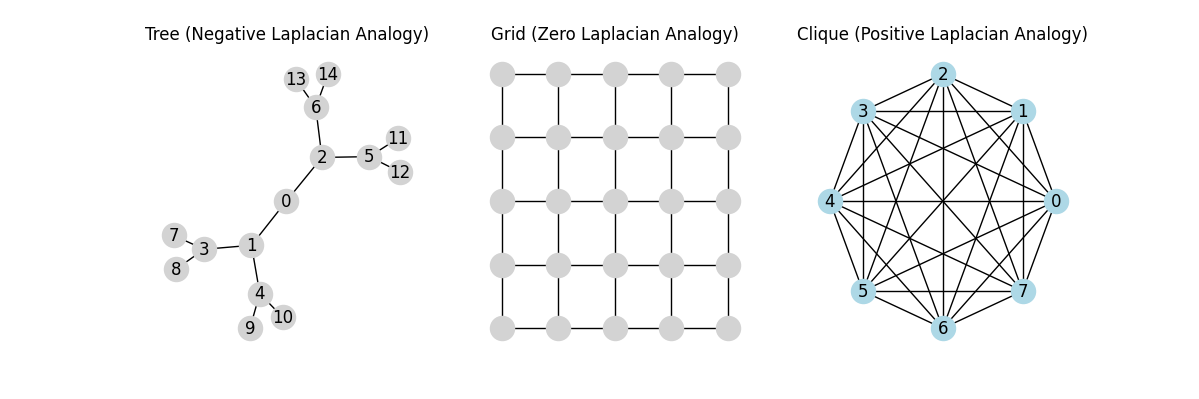}
    \caption{Shown are a tree graph exhibiting negative curvature except at leaf edges, an infinite graph with zero curvature throughout, and a complete graph with uniformly positive curvature respectively. It shows that Laplacian procure similar results as Ricci Flow.}
    \label{fig:laplacian graph}
\end{figure}
\begin{comment}
    \begin{figure}[h!]
    \hspace{-1cm}
       \centering
        \includegraphics[width=9cm, height=4cm]{illustrations/laplacian_manifold.png}
        \caption{Examples of Laplacian behavior on manifolds similar to Ricci flow; In the above figure, it shows manifolds with negative, zero and positive  respectively}
        \label{fig:laplacian manifold}
     
    \end{figure}
\end{comment}

\subsection{Edge Laplacian: Definition and Spectral Properties}

Let $G = (V,E)$ be an oriented graph with incidence matrix
$M \in \mathbb{R}^{|E| \times |V|}$. We define the (weighted) edge Laplacian
\begin{equation}
L := \operatorname{diag}(W_e)^{\frac{-1}{2}} \, M \, \operatorname{diag}(W_v) \, M^\top \operatorname{diag}(W_e)^{\frac{-1}{2}} ,
\label{eq:edge_laplacian}
\end{equation}
where $W_e \in \mathbb{R}^{|E|}$ and $W_v \in \mathbb{R}^{|V|}$ are diagonal matrices encoding edge and vertex weights, respectively. The operator $L$ acts on edge functions
$u \in \mathbb{R}^{|E|}$.

Throughout, we assume that
\[
W_e \succ 0,
\qquad
W_v \succ 0,
\]
i.e., both weight matrices are diagonal with strictly positive entries.

\subsubsection{Self-Adjointness}

Consider the symmetric edge Laplacian
\begin{equation}
L = \operatorname{diag}(W_e)^{-1/2} M\operatorname{diag}(W_v)M^\top \operatorname{diag}(W_e)^{-1/2}.
%\label{eq:edge_lap}
\end{equation}
Since symmetric normalization is used, the operator acts naturally on $\mathbb{R}^{|E|}$ endowed with the standard Euclidean inner product
\[
\langle x,y\rangle = x^\top y.
\]
\textbf{Proposition}:
The edge Laplacian $L$ is self-adjoint with respect to the standard Euclidean inner product.

\textit{Proof:}
Since $\operatorname{diag}(W_e)$ and $\operatorname{diag}(W_v)$ are diagonal matrices with positive entries implies they are symmetric. Therefore
\begin{align*}
L^\top &= \Big( \operatorname{diag}(W_e)^{-1/2} M\operatorname{diag}(W_v) M^\top \operatorname{diag}(W_e)^{-1/2} \Big)^\top \\ &= \operatorname{diag}(W_e)^{-1/2} M\operatorname{diag}(W_v) M^\top \operatorname{diag}(W_e)^{-1/2} \\ &=L.
\end{align*}
Hence $L$ is symmetric and therefore self-adjoint in finite dimensional vectorspace. It can also be seen using the direct inner product proof:
\[
\langle Lx,y\rangle = (Lx)^\top y= x^\top L^\top y=x^\top L y= \langle x, Ly\rangle.
\]

\subsubsection{Positive Semidefiniteness}\mbox{}

\noindent \textbf{Proposition}: The edge Laplacian $L$ is positive semidefinite.

\textit{Proof.}
Let $u\in \mathbb{R}^{|E|}$. Then $u^\top L u =$
\begin{align*}
u^\top \operatorname{diag}(W_e)^{-1/2} M\operatorname{diag}(W_v) M^\top \operatorname{diag (W_e)}^{-1/2}u.
\end{align*}
Define $z=M^\top\operatorname{diag}(W_e)^{-1/2}u.$ Then $u^\top L u= z^\top \operatorname{diag}(W_v)z.$ Since $W_v(i)>0$ for every vertex,
\[
z^\top \operatorname{diag}(W_v)z=\sum_{i \in [0,|V|] }W_v(i)\,z_i^2\ge 0.\]
Therefore, $u^\top L u \ge 0$,
and hence $L$ is positive semidefinite.

\subsubsection{Kernel and Structural Interpretation}\mbox{}

\noindent\textbf{Proposition.}
The kernel of $L$ is given by
\[
\ker(L)
=
\operatorname{diag}(W_e)^{1/2}\ker(M^\top).
\]
\textit{Proof:}
Since $L$ is positive semidefinite,
\[
u^\top L u =0
\]
if and only if $z= M^\top \operatorname{diag}(W_e)^{-1/2}u=0.$ Hence, $M^\top \operatorname{diag}(W_e)^{-1/2}u =0,$ which implies $\operatorname{diag}(W_e)^{-1/2}u \in \ker(M^\top).$ Multiplying by $\operatorname{diag}(W_e)^{1/2}$ yields
\[u \in \operatorname{diag}(W_e)^{1/2} \ker(M^\top).
\]
Therefore, $\ker(L)= \operatorname{diag}(W_e)^{1/2} \ker(M^\top).$ The kernel consists of weighted harmonic edge flows. Since $\operatorname{diag}(W_e)^{1/2}$ is invertible as all edge weights are strictly positive, $\ker(L)$ is isomorphic to $\ker(M^\top)$ and therefore retains the interpretation of divergence-free edge flows appearing in discrete Hodge theory.

\subsection{Intrinsic Stability of the Edge Heat Equation}

Since $L$ is self-adjoint and positive semidefinite, its spectrum is real and non-negative:
\[
0=\lambda_1 \le \lambda_2 \le \cdots \le \lambda_{|E|}.
\]

Consider the edge heat equation
\begin{equation}
\frac{\partial u}{\partial t} = -Lu.
\label{eq:edge_heat}
\end{equation}
\textbf{Theorem:}
Let $L$ be the edge Laplacian defined in \eqref{eq:edge_lap}. Then the operator $-L$ generates a strongly continuous contraction semigroup on $\ell_h^2(E)$. Furthermore, for every solution of \eqref{eq:edge_heat} lying in the orthogonal complement of $\ker(L)$, \[\|u(t)\|_2\le e^{-\lambda_2 t}\, \|u(0)\|_2, \] where $\lambda_2$ denotes the smallest positive eigenvalue of $L$.
\textit{Proof.}
Taking the Euclidean inner product of \eqref{eq:edge_heat} with $u(t)$ gives
\[\left\langle \frac{\partial u}{\partial t},u \right\rangle = -\langle Lu,u\rangle.\]
Since $ \left\langle \frac{\partial u}{\partial t}, u \right\rangle=\frac12\frac{d}{dt}\|u(t)\|_2^2,$ we obtain $\frac12\frac{d}{dt}\|u(t)\|_2^2 = -\langle Lu,u\rangle$. Restricting to the orthogonal complement of $\ker(L)$, the Rayleigh quotient characterization of $\lambda_2$ yields
\[ \langle Lu,u\rangle \ge \lambda_2\|u\|_2^2.\]
Hence, $\frac{d}{dt}\|u(t)\|_2^2 \le -2\lambda_2\|u(t)\|_2^2.$ Applying Gr\"onwall's inequality gives
\[ \|u(t)\|_2^2\le e^{-2\lambda_2 t} \|u(0)\|_2^2. \]
Taking square roots,
\[ \|u(t)\|_2\le e^{-\lambda_2 t}\|u(0)\|_2.\]
Therefore the edge heat flow is exponentially stable on $\ker(L)^\perp$. Consequently, for sufficiently small perturbations, the spectral gap persists, and the stability estimate for \eqref{eq:edge_heat} remains valid.

\textbf{Note:} 
\begin{itemize}
\item The proposed Laplace operator can be interpreted as a weighted and normalized lower Hodge1-Laplacian\cite{Lim2020Hodge}. In the absence of higher dimensional simplices, the Hodge 1-Laplacian reduces to $d_0d_0^*$ or ($-grad(div)$), which, under the edge by node incidence convection $M\in R^{|E|*|V|}$, gives $MM^T$. Our operator introduces vertex and edge weights and normalization through $W_v$ and $W_e$.

Moreover, we initially defined the Beltrami operator as $-div(grad)$, which only operates on node features on a Riemannian manifold which is similar to the Hodge 0-Laplacian while the Hodge 1-Laplacian acts on edge like quantities. But for graph without higher order simplices, reduces Hodge 1-Laplacian to a lower component $-grad(div)$. This motivates us to use the weighted and normalized edge-space Hodge Laplacian. 

\item Although discrete Ricci curvature notions may correlate with spectral gaps in certain graph families, the above stability result holds independently of any curvature interpretation. In contrast to discrete Ricci curvature approaches, no explicit curvature lower bound is imposed or required to ensure stability. While Ricci curvature can be interpreted as a geometric proxy for contraction properties of transport or diffusion processes, the edge-Laplacian framework yields stability as an intrinsic algebraic consequence of positivity and self-adjointness. In this sense, stability emerges naturally from the operator structure and is not enforced through artificial curvature assumptions.
    
\end{itemize}

\begin{algorithm}[h!]
\caption{Perturbation-Based Edge Reweighting and ARI Evaluation}
\label{alg:perturb-ari}
\textbf{Input:} Graph $G=(V,E)$, true labels $y$, number of iterations $T$, 
cut thresholds $\{c_1,\dots,c_K\}$ where $K$ denotes maximum checks to be done to obtain the optimum result, and perturbation step for initialisation $\epsilon$.\\
\textbf{Output:} ARI matrix $\mathrm{Ari} \in \mathbb{R}^{T \times K}$.

\begin{algorithmic}[1]

\State Compute node weights $W_v$ and incidence matrix $M \in \mathbb{R}^{|E|\times|V|}$.
\State Initialize edge  inverse weights

$W_e \leftarrow \mathbf{1}_{|E|} \ /\ (total\ no.\ of\ edges)$.
\State Compute edge Laplacian $L$,
 $L=diag(W_e)^{-1/2}M \ diag(W_v) \ M^T diag(W_e)^{-1/2}$
\State Step size $h_{\mathrm{opt}}$ such that  $ 0<h< \frac{2}{\lambda_{max}(L)}$
\State Compute perturbation Initialisation $u_0 \leftarrow \epsilon W_e$.

\For{$t = 1$ to $T$}
    
    \State Update perturbation: $u_h \leftarrow (I - hL_e)\,u_0$.
    \State $u_0=u_h$
    \State Update edge weights: $W_e \leftarrow W_e - u_h$.
    \State Normalize edge weights: $W_e \leftarrow W_e / \max(W_e)$.
    \State Assign updated $W_e$ to edges of $G$.

    \For{$j = 1$ to $K$}
        \State Remove all edges with weight $> c_j$ 
        \State New formed graph $H$ by the removal of edges.
        \State Compute predicted labels $\hat{y}$ by Louvain on $H$.
        \State Compute ARI: $\mathrm{Ari}[t,j] \leftarrow \mathrm{ARI}(y,\hat{y})$.
    \EndFor
    \State Storing the best cut to each iteration with best ARI for predicting the results.
\EndFor

\end{algorithmic}
\label{propAlg}
\end{algorithm}

%\newpage
\section{Results}
\subsection{Real World Dataset}

\subsubsection{Karate club Network\cite{zachary1977karate}}This dataset was collected by Wayne Zachary in 1970 from members of a University Karate Club. In this undirected graph nodes represent the members of the club and edges is the connection between two members. It has two classes which was formed after a conflict between the members.\\
\subsubsection{Football Network\cite{GirvanNewman2002} }This network is from the schedule of Division I games at the season fall 2000. Girvan and Newman has utilised this network for the community detection task.  This network can be partitioned into 12 conferences where each node represents the team and edge represents the match payed between the two teams. The idea for the communities in the dataset was that the match held in a conference(community) was more than the matches held between the teams in different conferences.
\subsubsection{Political blogs network\cite{10.1145/1134271.1134277} } This dataset was made during 2004 US Presidential elections considering the blogs posted by the liberals or conservative bloggers (represented as nodes). In this dataset, the edges were created if any blog is cited by the other blog.
\subsubsection{Email-EU-core network\cite{nr}} This network comprise of members of large European Research institute as nodes and the contact between two members through email as edges. Here, each individual member belongs to exactly one of the 42 departments in the institution.

\begin{table}[h!]
\centering
\caption{Real World Datasets}
\setlength{\tabcolsep}{4pt}
\begin{tabular}{|c|c|c|c|c|c|c|}
\hline
\textbf{DataSet} & \textbf{Vert.s} & \textbf{Edges} & \textbf{Class} & \textbf{Av. Deg.} & \textbf{CC} & \textbf{Diam.} \\
\hline
Karate Club & 34  & 78 & 2  & 4.5882 & 0.5706  & 5  \\
Football & 115 & 613  & 12 & 10.6608  & 0.4032  & 4  \\
Polblog & 1490 & 19025 & 2 & 27.3552 & 0.3202 & 8 \\
Email-EU core &  1005 & 16064 & 42 & 31.968 & 0.450 & 7 \\
\hline
\end{tabular}
\label{tab:dataset}
\end{table}

%\vspace{2cm}
\subsection{Experimental Results}
Applying our algorithm to the above datasets (Table \ref{tab:dataset}) yields good result. In order to evaluate the accuracy of the results, we used ARI (Adjusted Rand Index) as a metric to do it. It is the ratio of accuracy of the results obtained to the ground truth provided. 
\vspace{0.4cm}
\[\mathrm{ARI}=\frac{\sum_{i,j} \binom{n_{ij}}{2}-\frac{\sum_i \binom{a_i}{2}\sum_j \binom{b_j}{2}}{\binom{n}{2}}}{\frac{1}{2}\left(\sum_i \binom{a_i}{2}+\sum_j \binom{b_j}{2}\right)-\frac{\sum_i \binom{a_i}{2}\sum_j \binom{b_j}{2}}{\binom{n}{2}}}\]

where $n$ is the total number of data points, $n_{ij}$ denotes the number of data points common to cluster $i$ in the first partition and cluster $j$ in the second partition, $a_i=\sum_j n_{ij}$ and $b_j=\sum_i n_{ij}$ are the sizes of cluster $i$ and cluster $j$ in the first and second partitions, respectively.
%\vspace{0.5cm}

\begin{figure}[h!]
\centering

\includegraphics[width=\columnwidth,height=4cm]{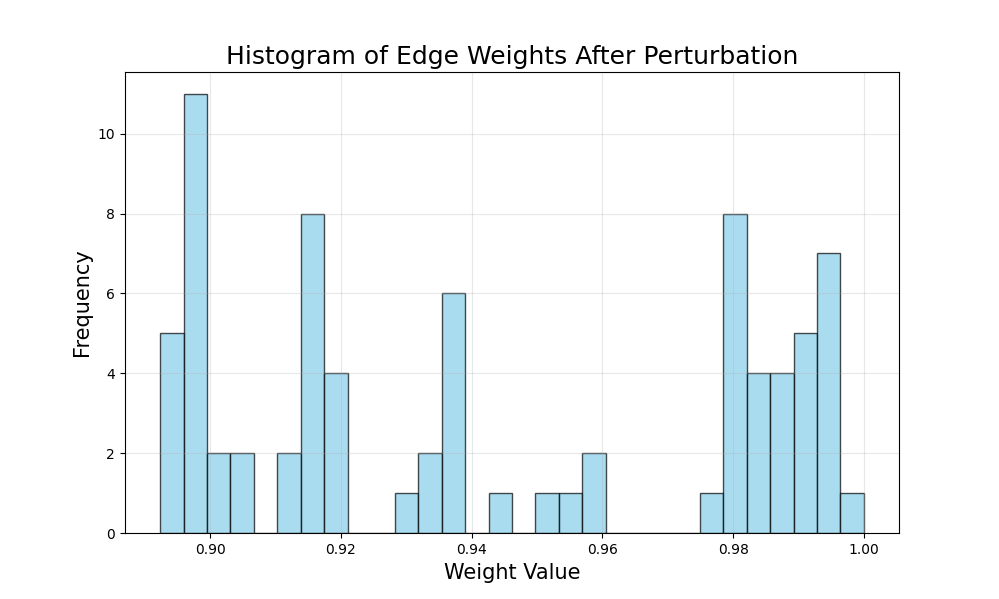}
\par\small (a) First iteration

%\vspace{0.1cm}
\includegraphics[width=\columnwidth,height=4cm]{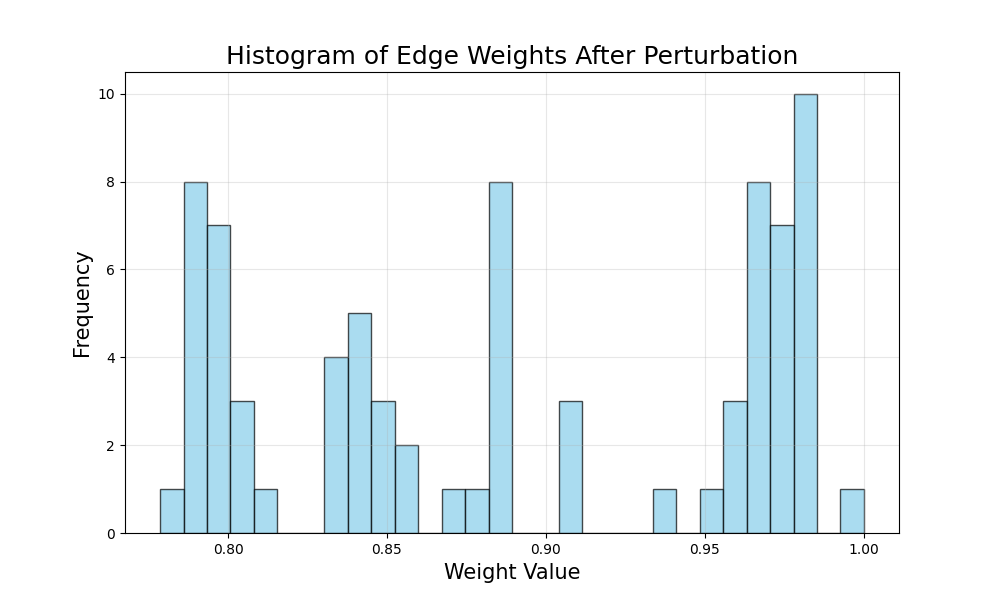}
\par\small (b) Second iteration

%\vspace{0.1cm}
\includegraphics[width=\columnwidth,height=4cm]{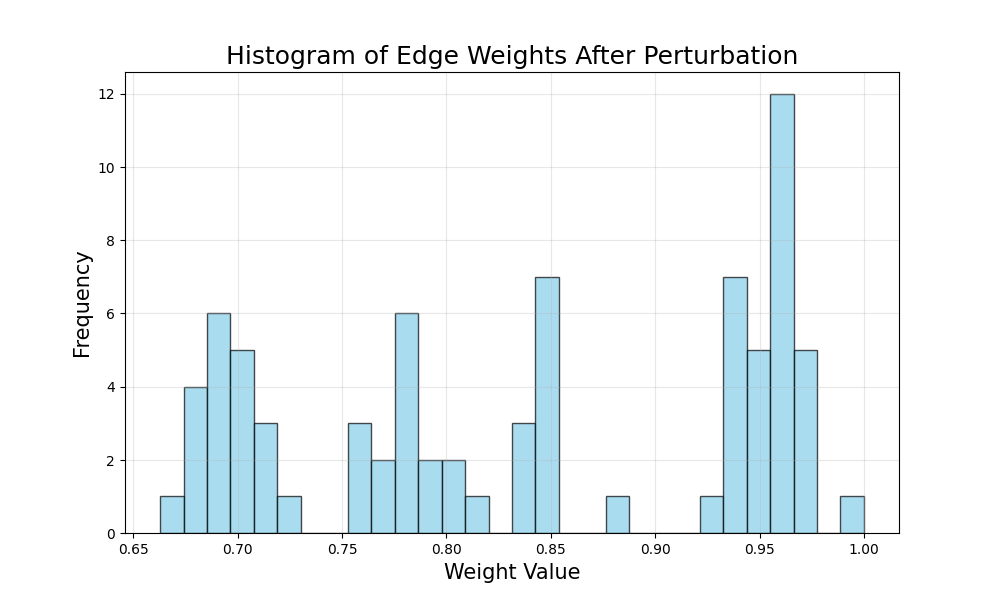}
\par\small (c) Third iteration

%\vspace{0.1cm}
\includegraphics[width=\columnwidth,height=4cm]{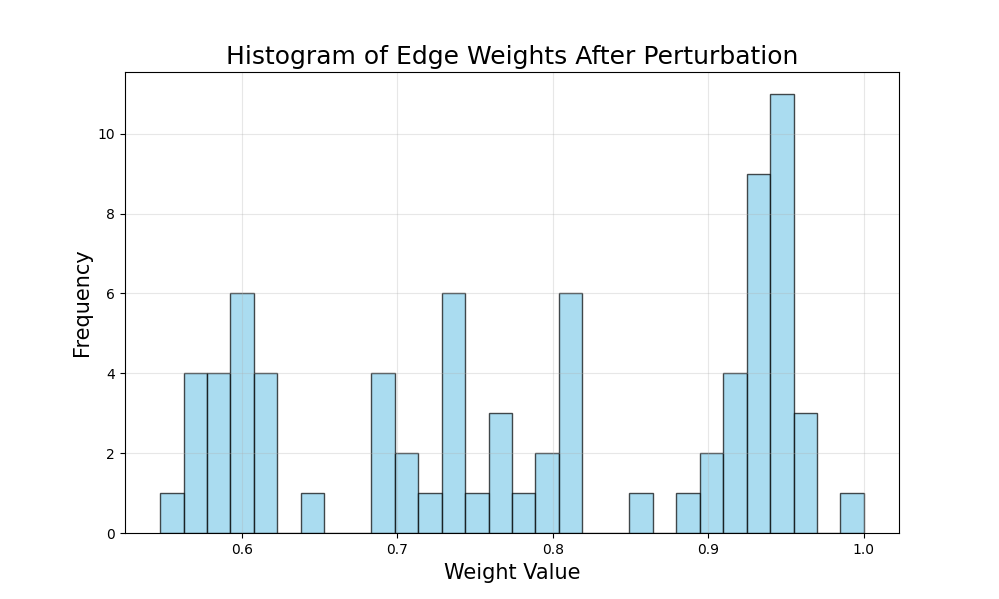}
\par\small (d) Fourth iteration

%\vspace{0.1cm}
\includegraphics[width=\columnwidth,height=4cm]{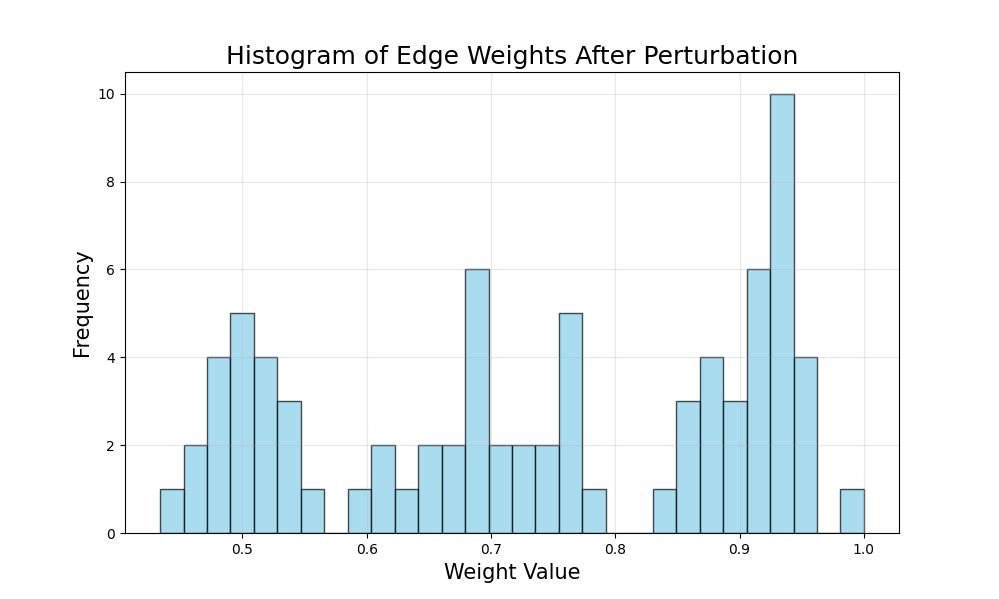}
\par\small (e) Fifth iteration

\caption{Histogram of edge-weight evolution across iterations for the Karate Club graph.}
\label{fig:histkarate}
\end{figure}

\begin{figure}[h!]
\centering
\begin{minipage}{0.5\textwidth}
    \centering
    \includegraphics[width=\textwidth, height=5cm]{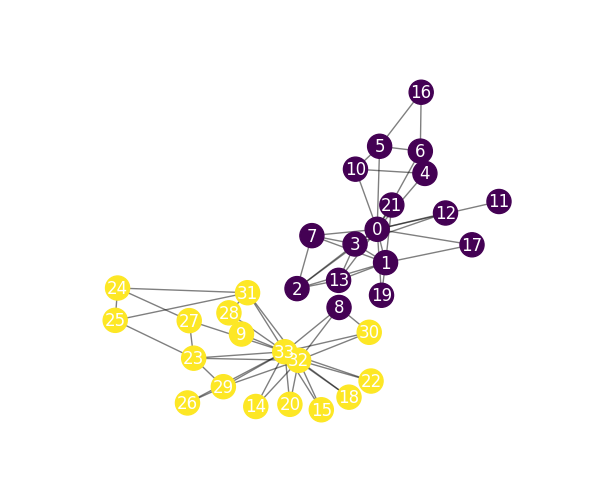}
    {(a) Karate Club Dataset}
\end{minipage}
\begin{minipage}{0.5\textwidth}
    \centering
    \includegraphics[width=\textwidth, height=5cm]{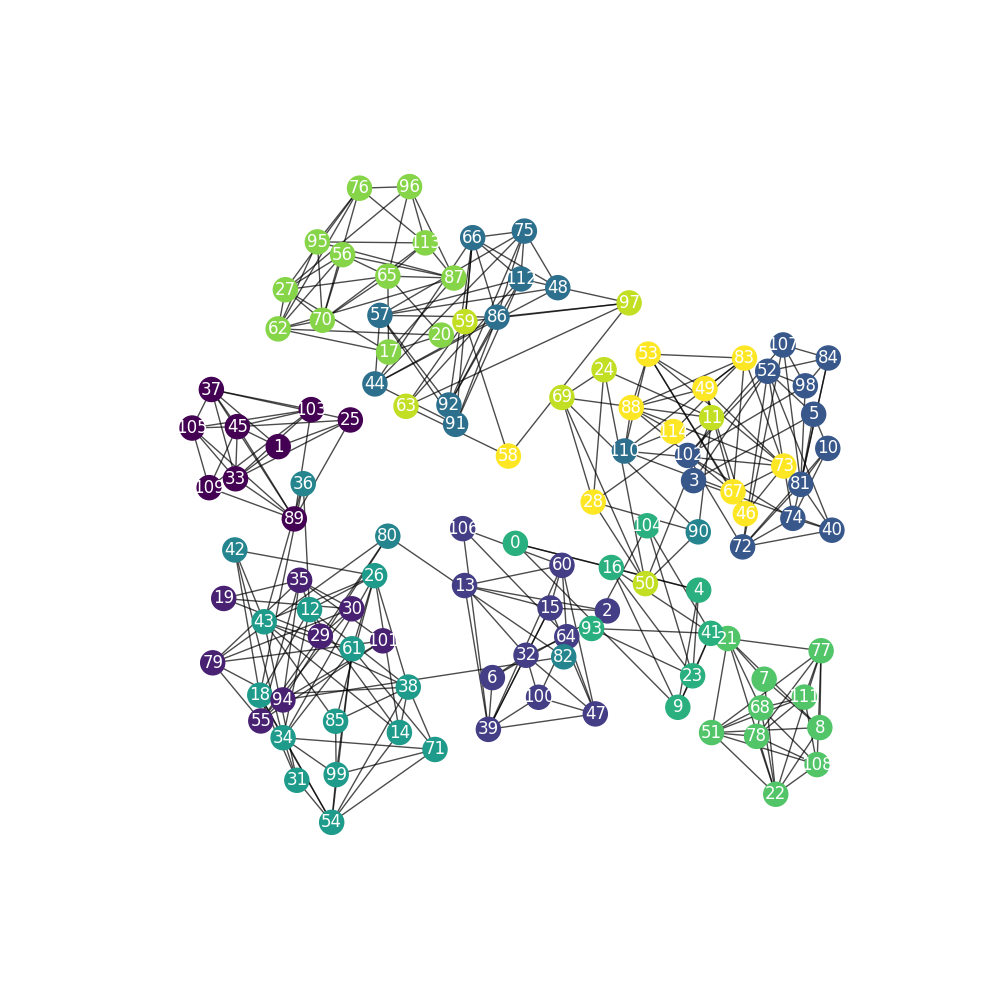}
    {(b) Football Dataset}
\end{minipage}

\caption{Representation of Clusters using the best ARIs obtained from the algorithm. Note: the graphs above are already divided into their original clusters and then the colored node represents the clustering done by the proposed algorithm of this paper.}
\label{fig:results}
\end{figure}

Within 3-5 iterations, we obtained good results, while discrete Ricci flow took more than 50 iteration to obtain the required results. We updated the weights at each iteration by removing the calculated perturbation and analyze the weight histogram (fig. \ref{fig:histkarate}) at each iteration to find the necessary cutoff. Finally, removing those edges with weights higher than the cutoff, we then applied Louvain algorithm \cite{blondel2008louvain} to find the community structure of the network. The cutoff was then obtained by considering the best ARI among the ARIs of different cutoffs (fig. \ref{fig:ARIVSCUT}). At the end, the results obtained are clearly visible in fig. \ref{fig:results}.  

For each dataset, the obtained results show that our method performs similar to Ricci flow, as illustrated in the figures. It also outperforms the approach of directly applying the Louvain algorithm, as well as Ricci flow for real world datasets. Except for the Polblog graph dataset, our algorithm produces good results for the other real world datasets. We performed analysis on synthetic dataset (Stochastic block model with 500 nodes and 2 communities) too, our Edge-laplacian based community detection method outperforms Ricci flow method and other traditional methods. This suggests that our method can be used for graph pruning or for computing graph embeddings in neural network–based dynamics. The main idea is to show that, instead of using a complex and computationally expensive algorithm, one can rely on a well defined Edge-Laplacian of a graph to capture curvature dynamics more efficiently.

\subsection{Computational efficiency}
\begin{figure}[h!]
\centering
\begin{minipage}{0.5\textwidth}
    \centering
    \includegraphics[width=\textwidth]{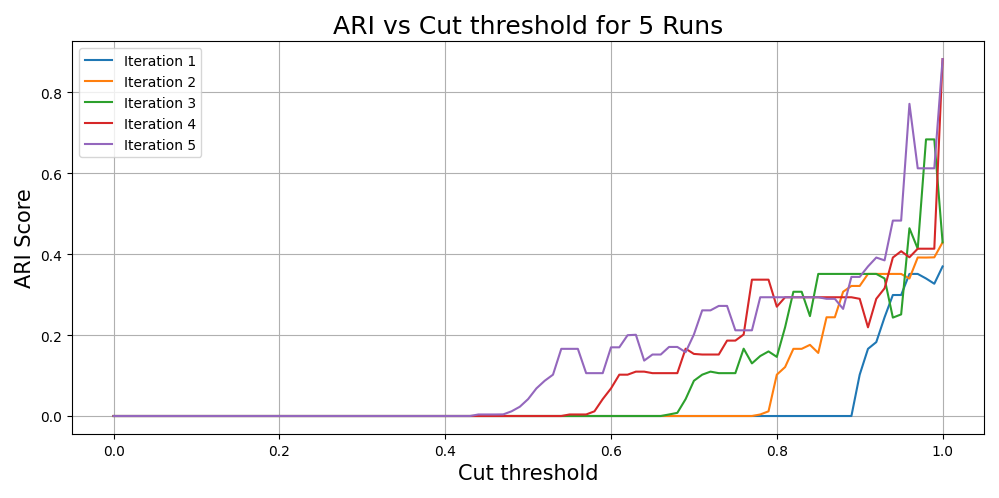}
    {(a) Karate Club Dataset}
\end{minipage}
\begin{minipage}{0.5\textwidth}
    \centering
    \includegraphics[width=\textwidth]{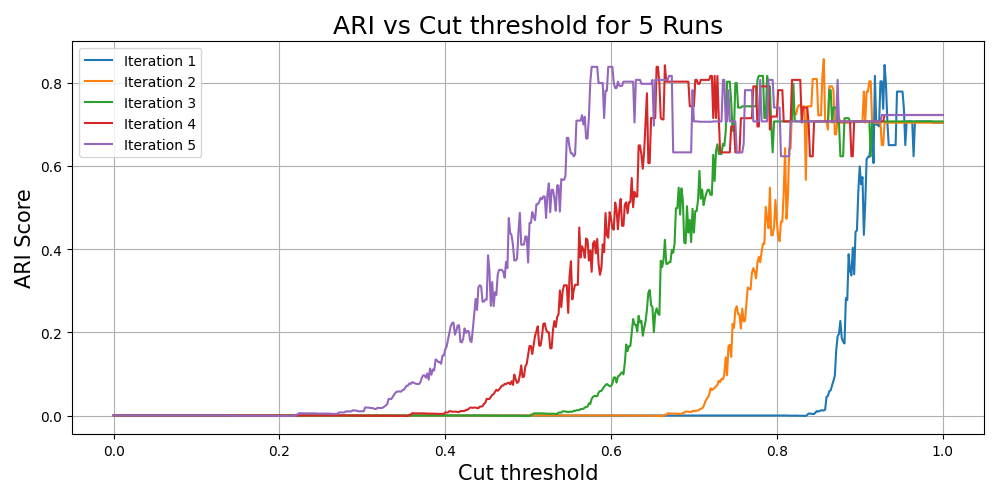}
    {(b) Football Dataset}
\end{minipage}
\begin{minipage}{0.5\textwidth}
    \centering
    \includegraphics[width=\textwidth]{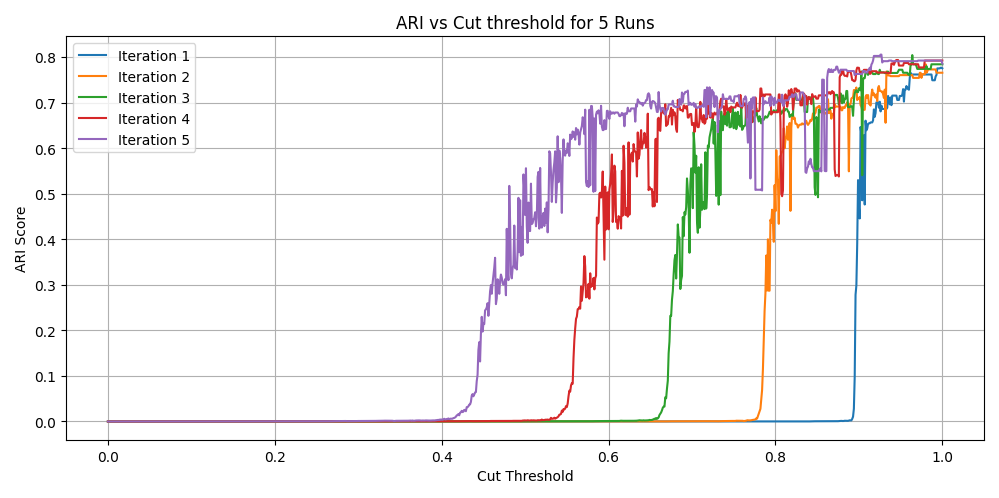}
    {(c) Polblog dataset}
\end{minipage}
\caption{Behavior of ARI for different cut off threshold for five iterations. Here we took 100, 600, 1000 equally discretized values in $[0,1]$ as cut values for the karate club graph, the football graph and polblog graph dataset respectively.}
\label{fig:ARIVSCUT}
\end{figure}

Ricci-flow-based methods relying on Ollivier-Ricci curvature require solving an optimal transport problem for every edge at each iteration, leading to a computational cost of $\mathcal{O}(m\bar d^3)$ per iteration, where $m$ is the number of edges and $\bar d$ is the average degree of two nodes whose Wasserstein distance (for worst case scenario it leads to solving a linear problem resulting in $\mathcal{O}(\bar d^3)$ computations) is to be calculated. In contrast, the proposed Edge-Laplacian-based flow is based solely on sparse matrix-vector multiplications and admits a linear $\mathcal{O}(m)$ complexity per iteration for calculating the edge Laplacian. While the dominant cost arises from evaluating K threshold values and performing Louvain clustering for each threshold, resulting in an overall complexity of approximately $\mathcal{O}(TKm)$ for $T$ flow iterations, where $K$ is treated as a small constant independent of graph size. As it can be clearly seen in fig. \ref{fig:ARIVSCUT} that the results are prominent after 0.5 cut off value which reduces the $K$ variable more. The proposed method scales approximately linearly with the number of edges, making it substantially more efficient than Ricci-flow methods based on optimal transport.
This reduction yields substantial speedups in practice.
\begin{table}[h!]
\centering
\caption{Comparison between Ricci flow and Edge based Laplacian}
\begin{tabular}{|c|c|}
\hline
\textbf{Ricci flow} & \textbf{Edge Based Laplace}\\
\hline
Loop over edges & Single local operator  \\
Compute curvature over edges & Apply global matrix \\
$m$ independent tasks & Sparse Matrix Vector multiplication\\
\hline
\end{tabular}
\end{table}

\begin{table}[h!]
\centering
\caption{ARI Results for Different Community Detection Methods}
\begin{tabular}{|c|c|c|c|}
\hline
\textbf{Dataset} & \textbf{Louvain} & \textbf{Ricci Flow} & \textbf{Proposed Flow} \\
\hline
Karate Club & 0.5089  & 0.7716  & \textbf{0.8823}  \\
Football  & 0.7041 &  0.8465 & \textbf{0.8569} \\
Polblog & 0.3768 & 0.7200 &  \textbf{0.8098}\\
Email-EU-core & 0.3260 & \textbf{0.4600} &  0.4306\\
\hline
\end{tabular}
\end{table}

\begin{figure}[h!]

    \centering
    \includegraphics[width=9cm, height=5cm]{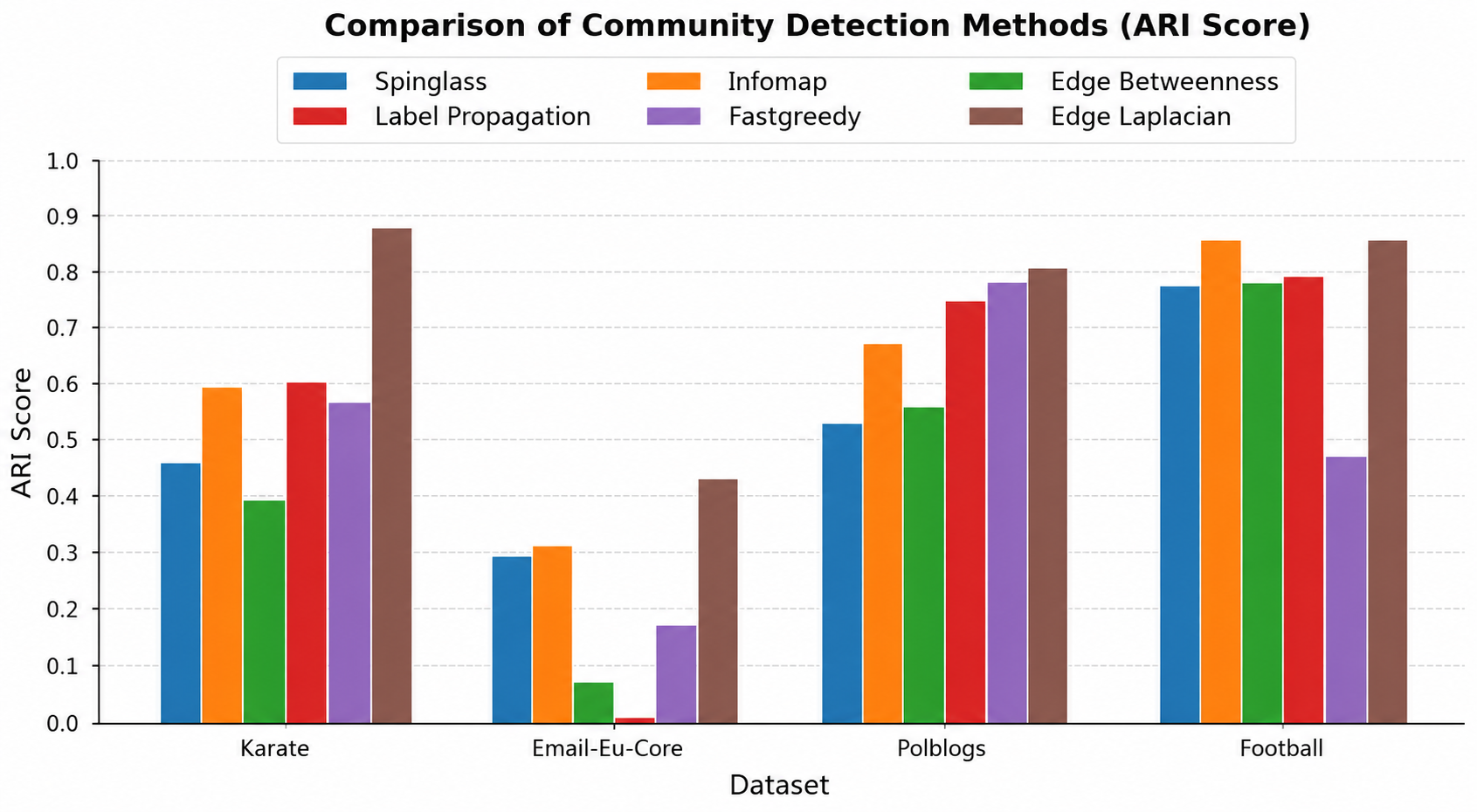}
    \caption{The accuracy of Edge laplacian for community detection on model networks measured by ARI.}
    \label{fig:ARIcomparison}
\end{figure}

\begin{figure}[h!]

    \centering
    \includegraphics[width=9cm, height=5cm]{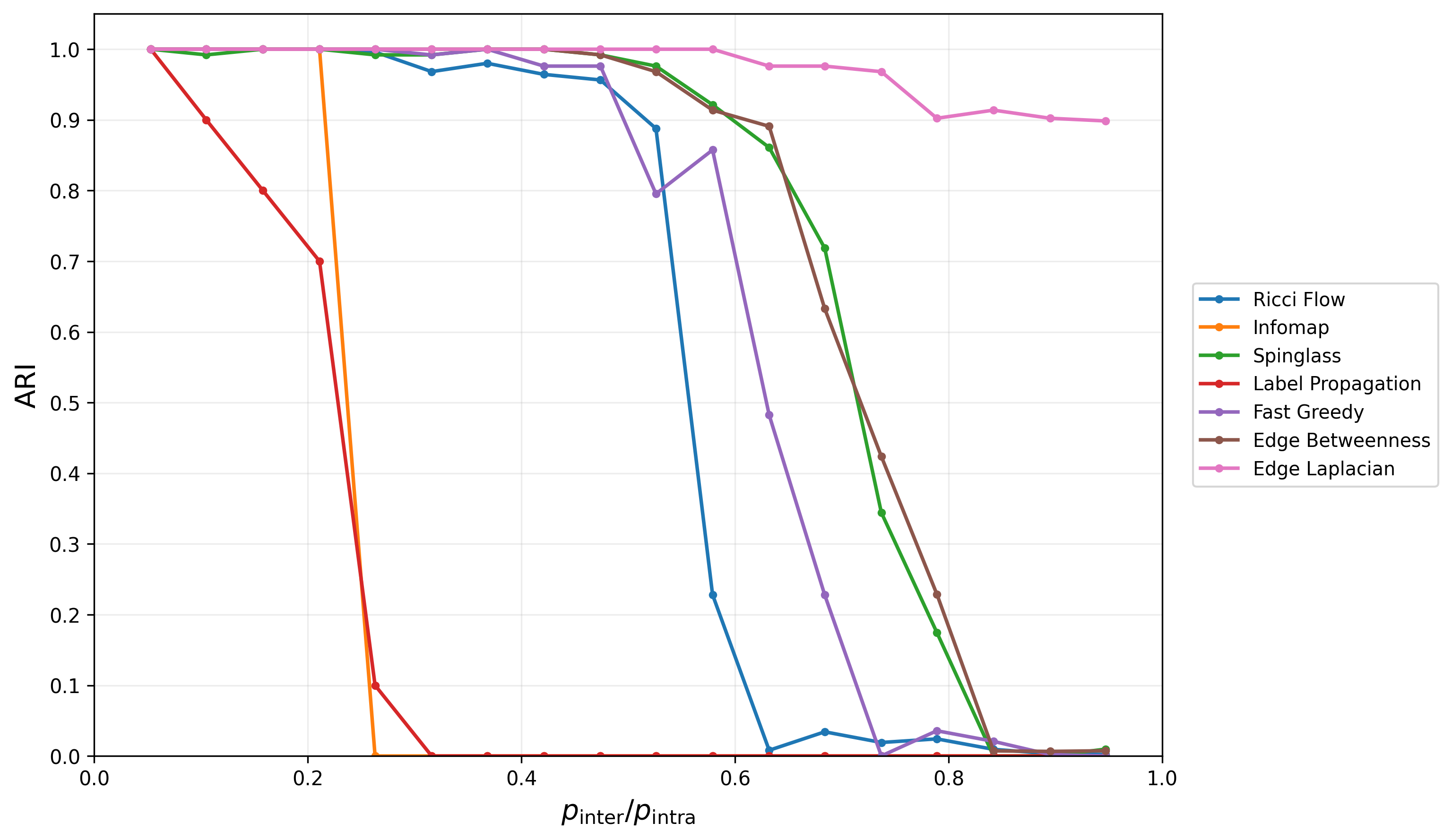}
    \caption{Considering the SBM graph of 500 nodes with two communities of equal size. The accuracy of Edge laplacian for community detection on SBM with various $p_{inter}/p_{intra}$ values, considering $p_{intra}=0.15$  is measured by ARI.}
    \label{fig:ARIcomparison}
\end{figure}


\vspace{1cm}

%\newpage
\section{Conclusion}
The results provided by the experiments state that the Laplacian can be chosen as an alternative to the discrete Ricci flow for certain graph datasets.
The results obtained are comparable to using a Ricci flow which can be computationally expensive considering the involvement of Wasserstein distance and the probability distribution associated with graph nodes or for evaluating the curvature for each edge involving square root computations over other edges on incident nodes. It may also be difficult to determine the underlying distribution of the nodes. In \cite{ni2019community}, a uniform distribution has been used but this may not work for all cases. Besides, Ricci curvature is known to produce spurious curvatures when there is none whereas the Edge- Laplacian is not an artificial computation. 
We observe that one can look to implement Edge-Laplacian method for community detection with cut threshold for a proper surgery instead of Louvian\cite{blondel2008louvain} to see the behavior of the graph. Finding a sustainable method of surgery can be considered as a future project. Moreover, the dataset considered here are only the small graphs, one can also look into large scale implementation of this method. It would be interesting to explore more into real world application considering our method.

\bibliographystyle{ieeetr}
\bibliography{myref}

\end{document}